\documentclass{fic-l}
                                                                                
\newtheorem{theorem}{Theorem}[section]
\newtheorem{lemma}[theorem]{Lemma}
\newtheorem{corollary}[theorem]{Corollary}
\newtheorem{proposition}[theorem]{Proposition}
                                                                                
\theoremstyle{definition}

\theoremstyle{remark}

\numberwithin{equation}{section}
                                                                                
\DeclareMathOperator{\Aut}{Aut}
\DeclareMathOperator{\Inn}{Inn}

\DeclareMathOperator{\Ref}{Ref}

\DeclareMathOperator{\odd}{odd}
\DeclareMathOperator{\eodd}{eodd}
\DeclareMathOperator{\Ang}{Ang}
\DeclareMathOperator{\FC}{FC}
\DeclareMathOperator{\GL}{GL}
\DeclareMathOperator{\Sym}{Sym}

\begin{document}
                                                                               
\title{The isomorphism problem for Coxeter groups}

\author{Bernhard M\"uhlherr}
\address{D\'epartement de Math\'ematiques\\
U.~L.~B. CP 216\\
Bd. du Triomphe\\
1050-Bruxelles\\
Belgium}
\email{bernhard.muhlherr@ulb.ac.be}

\subjclass{20F55; 51F55}

\begin{abstract}
By a recent result obtained by R. Howlett and the author
considerable progress has been made towards a complete solution of the 
isomorphism problem for Coxeter groups.
In this
paper we give a survey on the isomorphism
problem and  explain in particular  how the result mentioned above
reduces it to its `reflection preserving' version. Furthermore we
desrcibe recent developments concerning the solution of the latter.
\end{abstract}

\maketitle
                                                                                
\section{Introduction}

Coxeter groups are important in several mathematical
areas. It is therefore a bit surprising that the isomorphism
problem for those groups does not seem to have been considered
before the late 1990's. They only earlier reference known
to the author where this problem has been asked 
is \cite{Cohen}. The first major contributions to it 
are \cite{CD} and \cite{BMMN}. In \cite{CD} a rigidity
result is proved for a certain class of Coxeter groups. 
Rigidity  means that the Coxeter generating sets
are all conjugate. In \cite{BMMN} diagram twists
have been introduced. Those provide non-trivial
examples of non-rigid
Coxeter groups. The question about which Coxeter systems are rigid
arises naturally as well as the more general question about the 
isomorphism problem for Coxeter groups. 

The purpose of the present paper is to give a survey about what is
known at present about the isomorphism problem. The main motivation
for writing this survey is provided by a recent result obtained
by the author in collaboration with Bob Howlett. This result reduces
the isomorphism problem to its `reflection-preserving version'.
For the solution of the latter there is a conjecture
stated in \cite{BMMN}. Considerable progress towards a proof of
this conjecture was made in \cite{MW} and by recent work
of Pierre-Emmanuel Caprace in \cite{Ca05} there is reasonable hope
that this conjecture will be proved in the near future.
 Due to these facts
there is now a clear picture of what the solution of the
isomorphism problem should look like. In fact, at present there is 
a solution if one assumes that there are no irreducible spherical
residues of rank 3. We state two conjectures in Section \ref{sec5}.
The first one is known to be true for all Coxeter systems having
no $H_3$-subsystems; the second is a refinement of Conjecture 8.1
in \cite{BMMN} already mentioned.
Under the assumption that both conjectures are true,
we give an algorithm for the solution of the isomorphism problem.

\subsection*{Two versions of the isomorphism problem}

Let $W$ be a group and let $S \subseteq W$ be a set of involutions.
Then $M(S)$ denotes the square matrix $(o(ss'))_{s,s' \in S}$ where
$o(w)$ denotes the order of an element $w \in W$. The matrix $M(S)$
is called the {\it type} of $S$. 
As the elements of $S$ are involutions 
we have the following.

\begin{itemize}
\item[1.] For all $s,s' \in S$ we have $o(ss') \in {\bf N} \cup 
\{ \infty \}$;
\item[2.] for all $s \neq s' \in S$ we have $o(ss') = o(s's) \geq 2$;
\item[3.] for all $s \in S$ we have $o(ss) = 1$.
\end{itemize}

Hence, the matrix $M(S)$ is a symmetric square matrix 
with entries
in the set ${\bf N} \cup
\{ \infty \}$ where all entries on the main diagonal are equal to 
one and all remaining entries are strictly greater than one.
Such a matrix is called a {\it Coxeter matrix over $S$}.

Let $(W,S)$ be as above. We call $(W,S)$ a {\it Coxeter system}
(of type $M(S)$) if $\langle S \rangle = W$ and if the relations
$( (ss')^{o(ss')}=1_W)_{s,s' \in S}$ provide a presentation of $W$.
For a given Coxeter matrix $M = (M_{ij})_{i,j \in I}$ over a set $I$,
we define the {\it Coxeter group of 
type $M$} by setting 
 $W(M) := \langle I \mid ((ij)^{m_{ij}} =1 )_{i,j \in I} \rangle$.
It is a basic fact that the pair $(W(M),I)$ is a Coxeter system
of type $M$ (i.e. that $o(ij) = m_{ij}$ in $W(M)$ for all $i,j \in I$).

In this paper we will consider the isomorphism problem for 
finitely generated Coxeter groups. Thus, if we talk about
a Coxeter system $(W,S)$ or a Coxeter matrix $M$ over $I$ 
it is always understood that the sets $S$ and $I$ are finite.

Here are two versions of the isomorphism problem for Coxeter groups.

\smallskip
\noindent
{\bf Problem~1}
{\sl Given two Coxeter matrices $M$ and $M'$,
decide whether the groups $W(M)$ and  $W(M')$ are isomorphic.}

\smallskip
\noindent
{\bf Problem~2}
{\sl Given two Coxeter matrices $M$ and $M'$,
find all isomorphisms from $W(M)$ onto  $W(M')$.}

\smallskip
At first sight, Problem 1 seems to be a more natural question
than Problem 2. The latter is 
just a more general version of the first. 
Roughly speaking,
the solution of Problem 2 is equivalent to the solution of Problem 1
and a description of the automorphism group of $W(M)$ for any
Coxeter matrix $M$. This is in fact the main motivation to consider
Problem 2. It turns out that for certain 
Coxeter matrices $M$ a good understanding of 
the automorphism group of the 
group $W(M)$ is only possible
 if a solution of 
Problem 2 is available for all Coxeter matrices $M'$.

\subsection*{Content}

In Section \ref{sec2} we recall some definitions, fix notation
and mention some basic facts concerning Coxeter groups. In Section 
\ref{sec3} 
we will consider the rigidity problem for Coxeter groups.
This is an interesting special case of the isomorphism problem. In this
section we will provide examples of non-rigid Coxeter systems which
will play an important role later. Section \ref{sec4} is devoted to 
explaining
the results obtained in \cite{FHHM} and \cite{HM} and how these results
reduce the isomorphism problem to its `reflection-preserving version'
which will then be treated in Section \ref{sec5}. In Section 
\ref{sec6} we explain an
algorithm to solve the isomorphism problem under the assumption that
Conjectures 1 and 2 of Section \ref{sec5} hold. Finally, in Section 
\ref{sec7} we will make some remarks on the
automorphism groups of Coxeter groups.

\smallskip
\noindent
{\bf Remark:} It was mentioned above that there is no contribution to 
the isomorphism problem
for Coxeter groups before the late 1990's. Since then, however,
there are several publications concerning this subject. For instance,
Problem 1 has been solved completely in the case where $M$
is assumed to be even (i.e. no odd entries) by P. Bahls and M. Mihalik
(see \cite{Mi} and the references given there). 

In this
survey paper we do not attempt to give a systematic description 
of all contributions to the isomorphism problem for Coxeter
groups. We mention results (or consequences
of them)
whenever it will be  convenient. However, we try to include
all references on the subject in the bibliography. Thus, quite a few 
references will be mentioned only there.

\subsection*{Acknowledgement}

The content of this paper is based on my talk
at the
{\it Coxeter Legacy Conference} at Toronto in May 2004. I thank
the organizers for the invitation to 
present this survey at this conference.

\section{Preliminaries}

\label{sec2}

\subsection*{Coxeter diagrams}

With  a Coxeter matrix $M=(m_{ij})_{i,j \in I}$ we associate its
{\it diagram} $\Gamma(M)$. It is the edge-labelled graph $(I,E(M))$,
where the edge-set is 
$E(M) := \{ \{i,j \} \mid m_{ij} \geq 3 \}$ and where
an edge $ \{i,j \} \in E(M)$ has the label $ m_{ij}$. We do not
distinguish between a Coxeter matrix and its diagram since they 
carry the same information. We call a Coxeter matrix {\it irreducible}
if its associated Coxeter diagram is connected. An {\it irreducible 
component} of $M$ is a subset $J$ of $I$, which is a connected 
component of the diagram. A Coxeter matrix $M$ is called {\it spherical}
if $W(M)$ is finite. The irreducible spherical Coxeter diagrams
have been classified by H.S.M. Coxeter in \cite{Co}; we will use the 
Bourbaki notation for denoting them with the exception
that we denote rank 2 diagrams for the dihedral groups of order $2n$
by $I_2(n)$. Thus we have  the four series $A_n,C_n = B_n, D_n$
and $I_2(n)$ and the 6 exceptional diagrams $E_6, E_7, E_8, F_4, H_3$
and $H_4$.

An isomorphism from a Coxeter diagram $M=(m_{ij})_{i,j \in I}$
onto a Coxeter diagram $M'=(m'_{ij})_{i,j \in I'}$ is a graph 
isomorphism which preserves the edge-labels.

Let $M=(m_{ij})_{i,j \in I}$ be a Coxeter matrix over $I$
and let $J$ be a subset of $I$. Then we put 
$M_J := (m_{jk})_{j,k \in J}$ and $J^{\perp} := \{ k \in I \mid
m_{kj}=2 \mbox{ for all }j \in J \}$.

A Coxeter matrix $M$ is called {\it right-angled} if all edge-labels
of $\Gamma(M)$ are infinite; it is called {\it 2-spherical} if there
are no infinities; it is called {\it even} if there are no odd labels
 and it is called {\it of large type} if the diagram
is a complete graph (hence if there are no 2's in $M$).  

\subsection*{Coxeter systems}

Let $(W,S)$ be a Coxeter system. The set of its {\it reflections} is defined
to be the set $S^W:= \{ wsw^{-1} \mid s \in S \mbox{ and }w \in W \}$.
The {\it length} of $w \in W$ is the length of a shortest product
of elements in $S$ representing $w$; it is denoted by $l(w)$.
We call $(W,S)$ right-angled, 2-spherical, even  or of large type
if this is the case for $M(S)$.

We list some  facts about Coxeter systems which are important
in the sequel. Facts 1 and 2 are basic and can be found in
any standard reference on Coxeter groups (see \cite{NB} or \cite{JH});
Fact 3 is a non-trivial exercise in \cite{NB} but it follows also
from the fact that the Davis-complex of a Coxeter system is CAT(0);
Fact 4 is contained in \cite{Ri}; Fact 5 can be shown by considering
the geometric representation and Fact 6 is just an easy consequence
of the definition of a Coxeter system. 

\begin{itemize}
\item[1.] If  $J \subseteq S$, then 
$(\langle J \rangle, J)$ is a Coxeter system.

\item[2.] Let $J \subseteq S$ and $l:W \rightarrow {\bf N}$ be
the length function of $(W,S)$. Then the following are equivalent:
\begin{itemize}
\item[a)] $(\langle J \rangle, J)$ is finite;
\item[b)] there is an element $\rho_J$ such that $l(\rho_J) > l(x)$
for all $\rho_J \neq x \in \langle J \rangle$.
\end{itemize}
Moreover, if these two conditions are satisfied, then $\rho_J^2 = 1_W$.

\item[3.] If $X \leq W$ is a finite subgroup, then there exist 
$w \in W$ and $J \subseteq S$ such that $X^w \leq \langle J \rangle$
and such that $J$ is a spherical subset of $S$
(i.e. $\langle J \rangle$ finite).

\item[4.] Let $r \in W$ be an involution. Then there exist
$w \in W$ and $J \subseteq S$ such that
$J$ is spherical, $w \rho_J w^{-1} = r$ and such that
$\rho_J$ is central in $\langle J \rangle$.

\item[5.] Suppose that $J$ is a spherical subset of $S$ such that
$\rho_J$ is central in $\langle J \rangle$. Then the normalizer of
 $\langle J \rangle$ in $W$ and the centralizer of $\rho_J$
in $W$ coincide.

\item[6.] Let $(W,S)$ be a Coxeter system. Then each permutation
$\pi$ of $S$ which is an automorphism of $M(S)$ extends uniquely to
an automorphism $\gamma_{\pi}$ of $W$.

\end{itemize}

Let $(W,S)$ be a Coxeter system. By Fact 6 we can identify the stabilizer
of $S$ in $\Aut(W)$ with the group of automorphisms of $M(S)$;
this subgroup will be denoted by $\Gamma_S(W)$
and its elements are called the {\it graph-automorphisms} of $(W,S)$. 
The group $\Gamma_S$ has  trivial 
intersection with the group $\Inn(W)$ of inner automorphisms.
An automorphism of $W$ will be called {\it inner-by-graph} if it can be
written as a product of an inner automorphism and a graph-automorphism.

\section{Rigidity}

\label{sec3}

Let $G$ be a group and $R \subseteq G$ a set of involutions. Recall that
the  
Coxeter matrix $M(R)$ is called the type of $R$; the set $R$
is called {\it universal} if $(\langle R \rangle, R)$ is a Coxeter
system; it is called a {\it Coxeter generating set of $G$} if it is
universal and $G = \langle R \rangle$.

A Coxeter matrix $M$ is called {\it rigid} if for each Coxeter generating set
$R$ of $W(M)$ the Coxeter diagrams  $M(R)$  and $M$
are isomorphic. It is called {\it strongly rigid} if any two Coxeter generating
sets of $W(M)$ are conjugate in $W(M)$.

Clearly, strong rigidity implies rigidity. If a Coxeter 
diagram is (strongly) rigid, then we call the corresponding
Coxeter group and Coxeter system (strongly) rigid as well.

If one can show that the Coxeter diagram $M$ of Problem 1
is rigid, then this problem  is trivially solved. The answer is just that
the Coxeter diagram $M'$ has to be isomorphic to $M$.

Similarly, if one can show that the Coxeter diagram $M$ is
strongly rigid, then Problem 2 is solved. An isomorphism
onto $W(M')$ exists if and only if $M'$ and $M$
are isomorphic. Moreover, the automorphism group of
$W(M)$ is just the semi-direct product of
the group of inner automorphisms with the group of 
graph-automorphisms of $W(M)$; in other words: all automorphisms
of $W$ are inner-by-graph.

There are several interesting classes of Coxeter systems which
are not rigid. Before describing them we present some positive results.
The first is due to D. Radcliffe \cite{Ra}.

\begin{theorem} \label{thm1}
Right-angled Coxeter systems are rigid.
\end{theorem}

Although we fixed the convention that all Coxeter systems
in this paper are by definition of finite rank it is appropriate
to mention that the theorem above has been generalized to 
right-angled Coxeter systems of arbitrary rank by A. Castella
(see \cite{AC}).
The next result about strong rigidity is the
result of R. Charney and M. Davis already mentioned in the 
introduction (see \cite{CD}).

\begin{theorem} \label{thm2}
Let $(W,S)$ be a Coxeter system. If $W$ is capable of 
acting effectively,
properly and
cocompactly on some contractible manifold, then $(W,S)$ is strongly rigid. 
In particular, Coxeter groups of affine and compact hyperbolic type
are strongly rigid.
\end{theorem}

The next result is very recent. An important step towards a proof
of it was already made in \cite{HRT}; in the version presented here
it is a consequence of the main results in \cite{CM} and \cite{FHHM}.

\begin{theorem} \label{thm3}
Suppose that $(W,S)$ is irreducible, non-spherical and 2-spherical,
then $(W,S)$ is strongly rigid.
\end{theorem}

In the following we describe two ways to manipulate
the generating set of a given Coxeter system in order to produce
a new one whose type is possibly non-isomorphic to the type
of the original one. It is conjectured (and known to be true
in a lot of special cases) that Coxeter systems are rigid
up to these manipulations. 

\subsection*{Pseudo-Transpositions}

Let $k \geq 1$ be a natural number and put $n := 2(2k+1)$. We consider
the dihedral group $W$ 
 of order $2n$ as the group
of isometries preserving a regular $n$-gon in the euclidian plane.
Let $s,t \in W$ be two reflections whose axes intersect
in an angle $\frac{\pi}{n}$ and let $\rho$ be the central symmetry.
Then it is easily seen that $\{s,t \}$ and $\{s, tst, \rho \}$
are both Coxeter generating sets for $W$ of type $I_2(n)$ and
$I_2(2k+1) \times A_1$ respectively. Thus, the dihedral group
of order $2n$ is a non-rigid Coxeter group
because it has two Coxeter generating
sets of different types. This example is of course trivial
and a bit cheating because one of the two Coxeter matrices 
is not irreducible. However, it can be used to produce
more general examples by taking direct products or free products.
In \cite{HM} pseudo-transpositions have been introduced in order
to describe the general feature.

Let $(W,S)$ be a Coxeter system and let $\tau \in S$.
We call $\tau$ a {\it pseudo-transposition} if the following holds.

\begin{itemize}
\item[PT1] There is a unique $t \in S$ such that $o(\tau t) = 2(2k+1)$
for some natural number $k \geq 1$.
\item[PT2] For all $s \in S \setminus \{ \tau, t \}$ one has 
$o(\tau s) \in \{ 2, \infty \}$ and if $o(s \tau) = 2$, then
$o(st)=2$ as well.
\end{itemize}

The following is an easy observation about pseudo-transpositions.

\begin{lemma} \label{lem4}
Let $(W,S)$ be a Coxeter system, let $\tau \in S$ be a pseudo-transposition
of $(W,S)$ and let $t \in S$ be as in the definition above. Then
$S \setminus \{\tau \} \cup \{ \tau t \tau, \rho_{\{\tau,t \}} \}$
is a Coxeter generating set of $W$.
\end{lemma}

There is also another kind of pseudo-transpositions for Coxeter systems
based on the fact that the Coxeter groups $W(C_n)$ and 
$W(D_n \times A_1)$ are isomorphic for odd $n$. They yield also
non-isomorphic Coxeter generating sets in a similar way.
We refer to \cite{HM} for the details.

Let $(W,S)$ be a Coxeter system, let $\tau \in S$ be a pseudo-transposition
and let $R$ be the `new' 
Coxeter generating set as described in the lemma above.
Then we call the Coxeter system $(W,R)$ an
{\it elementary  reduction of $(W,S)$}. A Coxeter system $(W,S')$ will be
called a {\it reduction of $(W,S)$} if it can be obtained from $(W,S)$ 
by a sequence
of elementary reductions. Finally, we call $(W,S)$ {\it reduced}, if
there are no pseudo-transpositions. It is easy to see
that each Coxeter system has a reduced reduction.

Given a Coxeter diagram $M$ over a set $I$, then a Coxeter diagram
$M'$ over $I'$ is called an {\it elementary reduction of $M$} if there
is an elementary reduction of the Coxeter system $(W(M),I)$ whose
type is isomorphic to $M'$; we call $M'$ a {\it  reduction of $M$} if
$M'$ can be obtained from $M$ by a sequence of elementary reductions
and we call $M$ {\it reduced} if the system $(W(M),I)$ has no
pseudo-transpositions.

Clearly, any rigid Coxeter system has to be reduced in view of
Lemma \ref{lem4} above. The following result is due to M. Mihalik
\cite{Mi} and is based on earlier work of P. Bahls \cite{Ba0};
it states that the converse is true for even Coxeter systems.

\begin{theorem} \label{thm5}
An even Coxeter system is rigid if and only if there is no 
pseudo-transposition.
\end{theorem}

Note that this result generalizes Theorem \ref{thm1}.

\subsection*{Twistings}

In this subsection we describe twistings as they were introduced
in \cite{BMMN} and we give some further definitions concerning
them.

Let $(W,S)$ be a Coxeter system and let $J,K \subseteq S$. We call
the pair $(J,K)$ an {\it $S$-admissible} pair if the following holds.

\begin{itemize}
\item[AD1] $J$ is a spherical subset of $S$ and $K \cap (J \cup J^{\perp})
= \emptyset$.
\item[AD2] For all $k \in K$ and $l \in L:=S \setminus (J \cup J^{\perp}
\cup K)$ the order of $kl$ is infinite.
\end{itemize}

An $S$-admissible pair $(J,K)$ is called trivial if $K$ or $L$ are
empty. For a $S$-admissible pair $(J,K)$ we put 
$T_{(J,K)}(S) := J \cup J^{\perp}
\cup K \cup \{ \rho_J l \rho_J \mid l \in L \}$.

The following lemma is not too difficult to prove (see \cite{BMMN}).

\begin{lemma} \label{lem6}
Let $(W,S)$ be a Coxeter system and let $(J,K)$ be a
$S$-admissible pair. Then $T_{(J,K)}(S)$ is a Coxeter generating
set of $W$ which is contained in $S^W$.
\end{lemma}

Let $(W,S), (J,K)$ and $S':=T_{(J,K)}(S)$ be as in the previous
lemma. If $\rho_J$ is central in $\langle J \rangle$, then it is
easily verified that $M(S)$ is isomorphic to $M(S')$.
If $\rho_J$ is not central in $\langle J \rangle$, then
$M(S)$ is not isomorphic to $M(S')$ in the generic case. The following
example of such a situation was given in \cite{BM}.

\smallskip
\noindent
{\bf Example:}
Let $(W,S)$ be a Coxeter system such that $S= \{ s_1,s_2,s_3,s_4 \}$
and such that $o(s_1s_2) = o(s_2s_3) = o(s_3s_4) = 3$
and $o(s_1s_3)=o(s_1s_4)=o(s_2s_4) = \infty$. We put 
$J := \{ s_2, s_3 \}$ and $K := \{s_1 \}$. It follows
that $$S' := T_{(J,K)}(S) := \{ s_1' := s_1,
s_2' := s_2, s_3' := s_3, s_4' := s_2s_3s_2s_4s_2s_3s_2 \}$$
and that $o(s_1's_2') = o(s_2's_3') = o(s_2's_4') = 3$
and $o(s_1's_3')=o(s_1's_4')=o(s_3's_4') = \infty$. Thus
$M(S)$ and $M(S')$ are not isomorphic.

\smallskip
Let $S,R$ be Coxeter generating sets of a group $W$; we call
$R$ a {\it twist} of $S$ if there is a $S$-admissible pair $(J,K)$ such
that $R = T_{(J,K)}(S)$. It is readily verified that $R$ is a twist
of $S$ if and only if $S$ is a twist of $R$ and that $S^W = R^W$
in this case. 
A Coxeter generating set $S$ is called {\it twist-rigid} if there
are no non-trivial $S$-admissible pairs; i.e. if there are no
twists of $S$ which are not conjugate to $S$ in $W$.

Let $M$ be a Coxeter matrix over $I$. A Coxeter matrix $M'$
is called a {\it twist of $M$} if there is a twist $I'$ of $I$
in the Coxeter system $(W(M),I)$ such that $M(I')$ is isomorphic
with $M'$. As before one verifies that $M'$ is a twist of $M$
if and only if $M$ is a twist of $M'$.   
                                                                                 
We close this section with a result about strong rigidity for 
Coxeter groups. Obviously, if $(W,S)$ is a strongly rigid Coxeter system,
then $S$ has to be twist-rigid.
The following theorem provides the converse under the additional
assumption that all Coxeter generating sets $R$ of $W$ are
contained in $S^W$. In view of Corollary \ref{cor9} below there
are `a lot of examples' where this assumption holds.

\begin{theorem} \label{thm7}
Let $M$ be a non-spherical, irreducible Coxeter diagram over $I$ 
such that there is no subdiagram of type $H_3$. Suppose that
$I$ is a twist-rigid subset of $W(M)$ and that all
Coxeter generating sets of $W(M)$ are contained in $I^{W(M)}$. Then $M$
is strongly rigid.
\end{theorem}

This theorem was first proved in the large-type case ($m_{ij} > 2$
for all $i,j$) in \cite{MW}; 
 the result as it is stated above
has been obtained recently by P.-E. Caprace \cite{Ca05}.

\section{The reduction to the restricted isomorphism problem}
\label{sec4}

The restricted isomorphism problems for Coxeter groups are
the following:

\smallskip
\noindent
{\bf Problem 3:}
{\sl Given a Coxeter system $(W,S)$ and a Coxeter
matrix $M$, decide whether there is a Coxeter generating set
$R \subseteq S^W$ of $W$ such that $M(R) = M$.}

\smallskip
\noindent
{\bf Problem 4:}
{\sl Given a Coxeter system $(W,S)$ and a Coxeter
matrix $M$, find all Coxeter generating sets $R \subseteq S^W$
of $W$ with $M(R) = M$.}

\smallskip
In \cite{HM} Problems 1 and 2 of the introduction have been
reduced to Problems 3 and 4 respectively. This reduction
is based on the results on the 
{\it finite continuation of a reflection in a Coxeter group},
which have been obtained in \cite{FHHM}. The purpose of
this section is to describe the results obtained in both references.
The original motivation for the investigations in 
\cite{FHHM} was to find a tool to characterize reflections in
abstract Coxeter groups. We first  provide  
some examples, where an abstract Coxeter group does not determine
`its set of reflections'.

We have already seen examples, where an abstract Coxeter group
has different Coxeter generating sets yielding different sets of
reflections. If $(W,S)$ is not reduced and if $R$ is an elementary
reduction of $S$, then $S^W \not \subseteq R^W$ and 
$R^W \not \subseteq S^W$. We will now obtain further examples
by producing automorphisms of Coxeter groups which do not
preserve reflections. There are two kinds of such automorphisms, namely
$s$-transvections and $J$-local automorphisms. 

\subsection*{$s$-Transvections}
Let $(W,S)$ be a Coxeter system and let $s \in S$. We define the 
odd connected component  of $s$ 
in the diagram $\Gamma(S)$ to be the set of all elements 
$t \in S$ for which there is a path from $s$ to $t$ such that  
all its edge-labels are odd. We denote the odd component of
$s$ by $\odd(s)$ and we put $$\eodd(s) := \odd(s) \cup
\{ t \in S \mid o(tt') \neq \infty \mbox{ for some }t' \in \odd(s) \}.$$
Let $J_s$ denote the irreducible component of $\eodd(s)$ which
contains $s$ and let $K_s$ denote the union of 
all spherical irreducible components of $\eodd(s)$ which do
not contain $s$. 

Let $z$ be an element in the 
center of  $\langle K_s \rangle$. We define the mapping
$\theta_{s,z}: S \rightarrow W$ by setting $\theta_{s,z}(t)=tz$
if $t \in \odd(s)$ and by setting $\theta_{s,z}(t)=t$ for the 
remaining $t \in S$. One readily verifies that this mapping
extends to an involutory automorphism of $W$ and that
$sz$ is not contained in $S^W$. Hence $\theta_{s,z}(S)$ is 
a Coxeter generating set of $W$ providing a different set of
reflections.

The involutory automorphism described above is called an
{\it $s$-transvection} of the Coxeter system $(W,S)$. In fact, the
definition of an $s$-transvection given in \cite{HM}
is slightly more general. This is due to particular instances
which might arise when there are subsystems of type $C_3$. Due
to these instances the formal definition of
an $s$-transvection is somewhat involved and will be omitted here.
Nevertheless, we give an example of such a $C_3$-transvection because -
unlike for the other kinds of automorphisms - it is not
an `obvious automorphism easily seen from the diagram'.

\smallskip
\noindent
{\bf Example}
Let $(W,S)$ be a Coxeter system where $S = \{ s,t,t',c \}$
such that $o(st) = \\ o(st') =3$, $o(ct)=o(ct')=4$, $o(sc)=2$ and
$o(tt')= \infty$. Define $\theta: S \rightarrow W$
by setting $\theta(c) := c$, $\theta(s) := sc$,
$\theta(t) := stcsts$ and $\theta(t') := st'cst's$. One
verifies that 
$\theta$ extends uniquely to an involutory automorphism of $W$.

\subsection*{$J$-local automorphisms}
Let $(W,S)$ be a Coxeter system. A subset $J$ 
of $S$ is called a {\it graph factor}
of $(W,S)$ if $J$ is spherical and if for all
$t \in S \setminus J$ either 
$tj = jt$ for all $j \in J$ or $o(tj) = \infty$ for 
all $j \in J$. 

Let $J$ be a graph factor of $(W,S)$ and let $\alpha$ be an automorphism
of $\langle J \rangle$. Then it is readily verified that
there is a unique automorphism of $W$ stabilizing 
the subgroup $\langle J \rangle$, inducing $\alpha$ on it
and inducing the identity on $S \setminus J$. We call such an automorphism
a {\it  $J$-local automorphism}. 

This observation can be used to produce non-reflection preserving
automorphisms. There are lots of examples of
finite Coxeter groups, having automorphisms which
are not reflection preserving. Obvious examples
are the elementary abelian 2-groups.
A particularly interesting example is of course
the exceptional automorphism of $\Sym(6)$ which
is the Coxeter group of type $A_5$.

\subsection*{The finite continuation of a reflection}

Let $(W,S)$ be a Coxeter system. As $S$ is supposed
to be finite and as each finite subgroup of $W$ is conjugate
to a subgroup of some spherical standard parabolic subgroup it
follows that there is an upper bound for the order of any 
finite subgroup of $W$. This implies that there is for any subgroup
$X$ of $W$ a unique maximal normal finite subgroup of $X$ which we
denote by $O_{\text{fin}}(X)$.

Let $r \in W$ be an involution of $W$; by the result of Richardson
mentioned in Section \ref{sec2} (Fact 4) we know that $r$ is conjugate to
some $\rho_J$ for some spherical subset $J$ of $S$ and such that
$\rho_J$ is central in $\langle J \rangle$. Now one knows
that $N_W(\langle J \rangle) = C_W(\rho_J)$ (Fact 5) and hence
$\langle J \rangle$ is contained in $O_{\text{fin}}(C_W(\rho_J))$.
These considerations show that $r$ must be a reflection
if $O_{\text{fin}}(C_W(r)) = \langle r \rangle$. Hence we have found
a handy criterion which ensures that a given involution of an
abstract Coxeter group is a reflection for any Coxeter generating
set of that group.

This idea was the starting point for the results obtained
in \cite{FHHM}. It soon turned out that it is more convenient
to work with the {\it finite continuation $\FC(r)$} rather than
with the group $O_{\text{fin}}(C_W(r))$. This is defined to be 
the intersection of all maximal finite subgroups of $W$ containing
$r$. The main result of \cite{FHHM} is the following
theorem. Its proof
is based on a careful analysis of the centralizer of a reflection
which had been desrcibed in detail  in \cite{BB}.

\begin{theorem} \label{thm8}
Let $(W,S)$ be a Coxeter system
and let $s \in S$. Then $\FC(s)$ is known.  
Moreover, if  $\FC(s) =
\langle s \rangle$ , then $s$ is a reflection for each Coxeter 
generating set of $W$.
\end{theorem}

The description of $\FC(s)$ may become complicated if there
are subsystems of type $C_3$ or $D_4$. If this is not the
case, one can describe $\FC(s)$ by means of the subsets $J_s$
and $K_s$ defined in the paragraph on $s$-transvections as follows.

\begin{corollary} \label{cor9}
Let $(W,S)$ be a Coxeter system and suppose that there is no subsystem
of type $C_3$ or $D_4$. Let $s \in S$. If $J_s$ is spherical,
then $\FC(s) = \langle J_s \cup K_s \rangle$; in the remaining
cases one has $\FC(s) = \langle \{ s\} \cup K_s \rangle$.
In particular, if $K_s = \emptyset$ and $J_s$ is non-spherical,
then $s$ is a reflection for each Coxeter
generating set of $W$.
\end{corollary}

\subsection*{The reduction theorem}

Let $(W,S)$ be a Coxeter system. We call $s \in S$
{\it $\FC$-centered} if $\FC(s) = \langle J \rangle$ for some
$J \subseteq S$. A fundamental reflection
 might not be $\FC$-centered if there
are subsystems of type $C_3$ or $D_4$. Moreover,
the group of automorphisms of $W$ which stabilize the subset
$S^W$ is denoted by $\Ref_S(W)$.
We are now able to state the main result of \cite{HM}.

\begin{theorem}  \label{thm10}
Let $(W,S)$ be a reduced Coxeter system.
For each $\FC$-centered $s \in S$, let $T_s$ denote the group of
all $s$-transvections of $(W,S)$. For each
graph factor $J \subseteq S$ let $L_J$ denote the group
of all $J$-local automorphisms of $(W,S)$. Let $\Sigma$
be the subgroup of $\Aut(W)$ which is generated by all
$T_s$ and all $L_J$,
 where $s$ runs through
the $\FC$-centered
elements of $S$ and $J$ runs through the set of graph factors
of $(W,S)$.
Let $\tilde{\Sigma}$ be the subgroup
of $\Aut(W)$
which stabilizes $\FC(s)$
for all $s \in S$. Then we have the following:
\begin{itemize}
\item[a)] The group $\tilde{\Sigma}$ is finite and 
$\Sigma \leq \tilde{\Sigma}$.
In particular, $\Sigma$ is a finite subgroup of $\Aut(W)$.
\item[b)] Given a reduced Coxeter system $(W',S')$ and an isomorphism
$\alpha: W \rightarrow W'$, then there exists $\sigma \in \Sigma$
such that $\alpha(\sigma(S)) \subseteq S'^{W'}$.
\item[c)] The group $\Sigma$ (and hence also the group
$\tilde{\Sigma}$) is a finite supplement of $\Ref_S(W)$ in $\Aut(W)$.
\end{itemize}
\end{theorem}
 
Part b) of the theorem above says in particular, that if $(W,S)$
and $(W',S')$ are Coxeter systems which are both reduced and
if there is an isomorphism from $W$ onto $W'$, then there is
also an isomorphism between them which maps $S^W$
onto $S'^{W'}$. This yields the reduction of Problem 1
to Problem 3 for reduced Coxeter systems. 
Moreover, given any reduced Coxeter system
$(W,S)$, then its group of automorphism can be written as
$\Sigma \Ref_S(W)$, hence Problem 2 is reduced to Problem 4
for reduced Coxeter systems.
                                                  
\section{The restricted isomorphism problem}
\label{sec5}
In view of the reduction result described in the previous section
it suffices to solve Problems 3 and 4 in order to solve Problems
1 and 2 respectively. Thus we are led to the following question.

\smallskip
\noindent
{\bf Question:} 
{\sl Let $(W,S)$ be a Coxeter system and let $R \subseteq S^W$
be a Coxeter generating set of $W$. What can be said about $R$?}

\smallskip
We have to consider Coxeter generating sets whose elements
are reflections in a given Coxeter system. The following is
a first observation which can be shown by using the geometric
representation of a Coxeter group.

\begin{lemma} \label{lem11}
Let $(W,S)$ be a Coxeter system, let $R \subseteq S^W$ be a 
Coxeter generating set of $W$ and let $X \subseteq R$ be such that
$\langle X \rangle$ is finite. Then there exists a subset
$J$ of $S$ and an element $w \in W$ such that 
$\langle X \rangle^w = \langle J \rangle$. In particular,
if $r, r' \in R$ are such that $o(rr')=n \neq \infty$, then
there exist $s,s' \in S$ such that $o(ss') = n$ and 
such that the subgroups $\langle r, r' \rangle$ and 
$\langle s,s' \rangle$ are conjugate.
\end{lemma}

Let $(W,S)$ be a Coxeter system and let $R \subseteq S^W$ be
a Coxeter generating set. We call $R$ {\it sharp-angled}
with respect to $S$ if for any two reflections $r, r' \in R$
there exists $w \in W$ such that $\{r,r'\}^w \subseteq S$.

Let $W$ be the dihedral group of order $2n$ for some natural
number $n \geq 2$. We consider $W$ as the group of 
automorphisms of the regular $n$-gon in the euclidean plane.
Let $S= \{ s,t \}$, where $s$ and $t$ are reflections whose
axes intersect in an angle $\frac{\pi}{n}$. Given $r \neq r' \in S^W$,
then $\{ r,r' \}$ is sharp-angled with respect to $S$ if
the reflection axes of $r$ and $r'$ intersect in an angle
 $\frac{\pi}{n}$.

\subsection*{Angle-deformations}
Let $(W,S)$ be a Coxeter system, let $s \neq t \in S$
be such that $st$ has finite order, let
$x \in \langle s,t \rangle$ be such that $\langle s, xtx^{-1} 
\rangle = \langle s,t \rangle$ and put
$Y := S \setminus (\{ s,t \} \cup \{ s,t \}^{\perp})$.
Let $Y_s$ be the set of all $y \in Y$ for which there
exists a sequence $y_1,\ldots,y_k=y$ in $Y$ such that
$o(sy_1), o(y_1y_2), \ldots, o(y_{k-1}y_k)$ are finite and 
define $Y_t$ analogously. We define the mapping
$\delta_x: S \rightarrow W$ by setting $\delta_x(r) := r$
if $r \in S \setminus (Y_t \cup \{ t \})$ and $\delta_x(r) =
xrx^{-1}$ in the remaining cases.
The following is easy to verify.

\begin{lemma} \label{lem12}
If $Y_s \cap Y_t = \emptyset$ then $\delta_x$ extends uniquely
to an automorphism of $W$ which stabilizes the set $S^W$.
\end{lemma}

If $\{ s, xtx^{-1} \}$ is not sharp-angled with respect
to $\{ s,t \}$ and if $\delta_x$ is as above, then
$\delta_x(S)$ is not sharp-angled with respect to $S$. We therefore
call the automorphisms of the lemma above {\it angle-deformations}.  

The following result can be obtained by using rigidity of Fuchsian
Coxeter groups in a similar way as it was done in \cite{MW}.

\begin{proposition} \label{prop13} 
Let $(W,S)$ be a Coxeter system and suppose that
there is no 3-subset $J$ of $S$ such that $M(J) = H_3$. Let
$\Delta$ be the group generated by all angle deformations of
$(W,S)$. Given a Coxeter generating set $R \subseteq S^W$,
then there exists $\delta \in \Delta$ such that
$\delta(R)$ is sharp-angled with respect to $S$.
\end{proposition}

In view of the previous proposition the following conjecture
is known to be true for Coxeter systems having no subsystem
of type $H_3$.

\smallskip
\noindent
{\bf Conjecture 1:} 
{\sl Let $(W,S)$ be a Coxeter system and $R \subseteq S^W$
be a Coxeter generating set.
Then there exists an automorphism $\alpha$ of $W$ such that
$\alpha(S^W) = S^W$ and such that $\alpha(R)$ is sharp-angled 
with respect to $S$.}

\subsection*{Twist-equivalence}

Let $(W,S)$ be a Coxeter system and let $R \subseteq S^W$
be a Coxeter generating set of $W$. Recall that $R' \subseteq S^W$
is called a twist of $R$ if there is an $R$-admissible pair
$(J,K)$ such that $R' = T_{(J,K)}(R)$.  Moreover,
$R'$ is a twist of $R$ if and only if $R$ is a twist of $R'$.
By taking the transitive closure we obtain an equivalence relation
on the set of the Coxeter generating sets contained in $S^W$
which is called {\it twist-equivalence}.

If $R'$ is a twist of $R \subseteq S^W$, then $R' \subseteq R^W$
and $R'$ is sharp-angled with respect to $R$. Hence, if $R'$ is 
twist-equivalent with $R \subseteq S^W$, then $R' \subseteq R^W$
and $R'$ is sharp-angled with respect to $R$. There is some evidence
that the converse is also true. This is the content of the 
conjecture below. This conjecture is a refinement of Conjecture 8.1
in \cite{BMMN}.

\smallskip
\noindent
{\bf Conjecture 2:}
{\sl Let $(W,S)$ be a Coxeter system and
$R \subseteq S^W$ a Coxeter generating set of $W$ which is
sharp-angled with respect to $S$. Then $R$ is twist-equivalent
to $S$.}

\smallskip
At present, the following two theorems are known by recent work of 
P.-E. Caprace. The first improves earlier results obtained
in \cite{BMMN},
and \cite{MW}.

\begin{theorem} \label{thm14}
Conjecture 2 holds for all 
Coxeter systems which  do not contain an irreducible spherical
subsystem of rank 3.
\end{theorem}

\begin{theorem} \label{thm14a}
If $(W,S)$  is  a Coxeter system such that $M(J^{\perp})$
is 2-spherical for each spherical subset $J$ of $S$,
Conjecture 2 holds for $(W,S)$.
\end{theorem}

The main tool to prove Conjecture 2 in the references
above is known to the experts 
as `Kac Conjugacy Theorem for root bases'. This
theorem is proved in \cite{Kac} for affine and compact hyperbolic groups.
A proof for all Coxeter groups is  given in \cite{HRT}.

\section{The solution of Problem 1}
\label{sec6}

Let $M$ be a Coxeter diagram over a set $I$. Recall
that $M'$ is called a twist of $M$ if there 
is a twist $I'$ of $I \subseteq W(M)$ such that
$M(I')$ is isomorphic to $M'$. Again, $M'$ is a twist
of $M$ if and only if $M$ is a twist of $M'$ and by
taking the transitive closure we obtain an equivalence relation
on the set of Coxeter matrices which is  called 
twist-equivalence as well.

The following lemma is easy to prove.

\begin{lemma} \label{lem15}
Let $(W,S)$ be a Coxeter system and let $M$ be a Coxeter matrix.
Then the following are equivalent.

\begin{itemize}
\item[a)] There exists a Coxeter generating set $R \subseteq S^W$
such that $M(R)$ is isomorphic to $M$ and such that $R$
is twist-equivalent to $S$.
\item[b)] The matrices $M(S)$ and $M$ are twist-equivalent.
\end{itemize}
\end{lemma}

Using the previous lemma one obtains the following theorem, which
yields the solution of Problem 3.

\begin{theorem} \label{thm16}
Let $(W,S)$ and $(W',S')$ be Coxeter systems and suppose that
Conjectures 1 and 2 hold for $(W,S)$. Then the following
are equivalent.

\begin{itemize}
\item[a)] $M(S)$ and $M(S')$ are twist-equivalent.
\item[b)] There exists an isomorphism $\alpha:W' \rightarrow W$
such that $\alpha(S') \subseteq S^W$
\end{itemize}
\end{theorem}

We recall that a Coxeter system $(W,S)$ is reduced if the set $S$
contains 
no pseudo-transposition, that there is a natural
notion of a Coxeter system or a Coxeter matrix to be
a reduction of another and that it is always possible to
produce a reduced reduction of a Coxeter system or Coxeter matrix
by an easy algorithm. Now the previous theorem and 
Theorem \ref{thm10} yield the following.

\begin{theorem} \label{thm17}
Let $M$ and $M'$ be irreducible Coxeter matrices of rank at least 3 
and let 
$(W,S)$ be a Coxeter system of type $M$.
If Conjectures 1 and 2 hold for $(W,S)$, then 
                             the following are
equivalent.

\begin{itemize}
\item[a)] The groups $W(M)$ and $W(M')$ are isomorphic. 
\item[b)] If $M_1$ is a reduced reduction of $M$ and
if $M_1'$ is a reduced reduction of $M'$, then
$M_1$ and $M_1'$ are twist equivalent.
\end{itemize}
\end{theorem}

In view of Theorem \ref{thm14} and Proposition
\ref{prop13} we have the following corollary.
  
\begin{corollary} \label{cor18} 
Let $M$ and $M'$ be Coxeter matrices and suppose 
that $M$ has no subdiagram of type $A_3,C_3$ or $H_3$, 
then the following are equivalent: 
\begin{itemize} 
\item[a)] The groups $W(M)$ and $W(M')$ are isomorphic.
\item[b)] If $M_1$ is a reduced reduction of $M$ and
if $M_1'$ is a reduced reduction of $M'$, then
$M_1$ and $M_1'$ are twist equivalent.
\end{itemize}
\end{corollary} 

\section{On automorphisms of Coxeter groups}
\label{sec7}

The previous section shows that there is---under the
hypothesis that Conjectures 1 and 2 are true---a satisfactory
solution of Problem 1. Unfortunately, we cannot offer a satisfactory 
description of the automorphism groups of Coxeter groups under
the same assumptions which would yield 
a solution
of Problem 2 as well. In fact, the author has serious doubts
whether such a handy description exists in the general case.
Nevertheless there are several natural subgroups of the automorphism
group of a Coxeter group which are quite well understood.
In most of the `interesting' cases, the understanding of these subgroups
suffices to understand the group of automorphisms as a whole.
Our discussion will be restricted to those subgroups.
Before going more into the details we would like to mention
that the automorphism groups of Coxeter groups had been determined
in various special cases.

\begin{itemize}
\item[1.]  A presentation of the
 automorphism groups of right-angled Coxeter groups was given
in \cite{Mu98}. This work was based on the results obtained in
\cite{JT} and the latter is a far reaching generalization of 
the result in \cite{LJ}.
\item[2.] The automorphism groups of 2-spherical Coxeter groups
are `trivial' (i.e. all automorphisms are 
inner-by-graph) if there is no direct factor which is spherical.
This result was accomplished in \cite{CM} and \cite{FHHM}.
A `virtual' result in this direction 
has been obtained already in \cite{HRT} and
the main tool developed there was used again in \cite{CM}.
\item[3.] The automorphism groups of several classes of Coxeter groups
which are `almost spherical' have been described in \cite{WNF801},
\cite{WF}, \cite{FH1} and \cite{FH2}. In \cite{FH2} a complete
description of the automorphism groups of the irreducible spherical
Coxeter groups is given.  
\end{itemize}

Given an abstract Coxeter group $W$, then there is always a 
Coxeter generating set $S \subseteq W$ such that $(W,S)$
is reduced. Thus, there is no loss of generality if we consider
only reduced  Coxeter systems in this section.
Let $(W,S)$ be a reduced Coxeter system.
We define the following subgroups:

\begin{itemize}
\item[1.] $\Ref_S(W) := \{ \alpha \in \Aut(W) \mid \alpha(S^W) = S^W \}$,
\item[2.] $\Ang_S(W) := \{ \alpha \in \Ref_S(W) \mid \alpha(S) 
\mbox{ sharp-angled with respect to }S \}$,
\item[3.] $\tilde{\Sigma}_S(W) := \{ \alpha \in \Aut(W) \mid 
\alpha(\FC(s)) = \FC(s)
\mbox{ for all }s \in S \}$
\item[4.] $\Gamma_S(W) := \{ \alpha \in \Aut(W) \mid \alpha(S) = S \}$
\end{itemize}

In view of Theorem \ref{thm10} we have 
$\Aut(W) = \tilde{\Sigma}_S(W) \Ref_S(W)$ and the group
$\tilde{\Sigma}_S(W)$ is a finite group. Thus, there is 
a finite supplement of $\Ref_S(W)$ in $\Aut(W)$. There is 
the natural question about minimal supplements (or even complements)
of $\Ref_S(W)$ in $\Aut(W)$.
The example of the Coxeter group of type $A_1^k$ shows that
there are not always complements. However, a careful analysis
of several special cases provides some evidence for the following
conjecture.

\smallskip
\noindent
{\bf Conjecture 3:}
{\sl Let $(W,S)$ be a reduced Coxeter system. Then
there exists a subgroup $\Omega \leq \tilde{\Sigma}_S(W)$ such that
$\Pi:= \Omega \cap \Ref_S(W) \leq \Gamma_S(W)$ and such that
$\Omega$ is a supplement
of $\Ref_S(W)$ in $\Aut(W)$. Moreover, there is a normal 2-subgroup
$U$ of $\Omega$ and a complement of $L$ of $U$ in $\Omega$
such that $L=L_1 \times L_2 \times \ldots L_k$ where
$L_i$ is isomorphic to $\GL(n_i,2)$ for some natural number $n_i$
for $1 \leq i \leq k$ and $\Pi \cap L_i$ is just the set of
permutation matrices.}

There is a canonical candidate for the choice of the 
group $\Omega$ and based on this choice
the validity of the conjecture is not
difficult to see in several special cases.
However, the arguments become somewhat involved in the general
case. 

\subsection*{Reflection-preserving automorphisms}

As $\Ref_S(W)$ has a finite supplement, a big part of 
$\Aut(W)$ is understood if $\Ref_S(W)$ is understood. A
first observation is that $\Ang_S(W)$ is 
a normal subgroup of finite index in
$\Ref_S(W)$ and therefore a similar remark holds for $\Ang_S(W)$.
We do not know whether $\Ang_S(W)$ always has a finite supplement
in $\Ref_S(W)$ but we believe that there are examples where
this is not the case. If there is no $H_3$-subdiagram, then
the group $\Ref_S(W)$ is generated by the angle-deformations
of $(W,S)$ and $\Ang_S(W)$. We expect this to be true in general
with a suitable definition of angle-deformations in the case
where there are $H_3$-subdiagrams. 

In the following we will consider the group $\Ang_S(W)$.
Let $${\bf R} := \{ R \subseteq S^W \mid R \mbox{ sharp-angled 
Coxeter generating set of $W$ with
respect to }S \}$$ and call two elements $R \neq R'$ in ${\bf R}$
adjacent if one is a twist of the other. This yields 
a graph which we call ${\bf C}$. Conjecture 2 is equivalent
to the statement that the graph ${\bf C}$ is connected. 

We consider first the special case where $M(S)$ is even
in which case Conjectures 1 and 2 are known to be true.
If $M(S)$ is even, there is for each 
neighbor $R$ of $S$ in the graph ${\bf C}$ a canonical
involution $\theta_R$ in $\Ang_S(W)$ which switches $S$ and $R$.
Setting  
$X:= \langle \theta_R \mid R \mbox{ neighbor of }S \rangle$,
one verifies that ${\bf C}$ is the Cayley graph of $X$ with respect
to this generator set and that $\Gamma_S$ is a complement of 
$X$ in $\Ang_S(W)$. It is probably possible to generalize
the arguments given in \cite{Mu98} in order to give a presentation
of the group $\Ang_S(W)$. The key ingredient of such
a generalization would be the observation that the group
$\Ang_S(W)$ is something like a `generalized Coxeter group'
as it is in the right-angled case.

Let's consider the general case under the assumption that Conjecture 2
holds. 
The situation becomes more complicated.
The graph ${\bf C}$ is no longer the Cayley graph of a group
but of a groupoid. We do not go into the details here. But it
is worth mentioning that a similar situation occurs if
one is interested in the normalizer of a parabolic subgroup
in a Coxeter group. These normalizers had been described
in \cite{Bo} and \cite{BH} in a satisfactory way. 
The key observation in \cite{BH} is that they
are finite index subgroups of a groupoid which one might call
a Coxeter groupoid in view of its properties which are quite
similar to those of Coxeter groups. 
We believe,
that a presentation of $\Ang_S(W)$ can be given by using analogous 
ideas. It would be based on the observation that the graph
${\bf C}$ is the Cayley-graph of a generalized Coxeter groupoid
of which $\Ang_S(W)$ is a subgroup of finite index. 
However, a concrete description of such a presentation
might become rather involved.

\end{document}